\documentclass[12pt, reqno]{amsart}
\usepackage{amsmath, amstext, amsbsy, amssymb, amscd}

\newtheorem{theorem}{Theorem}[section]
\newtheorem{lemma}[theorem]{Lemma}

\newtheorem{corollary}[theorem]{Corollary}
\theoremstyle{definition}
\newtheorem{definition}[theorem]{Definition}

\theoremstyle{remark}
\newtheorem{remark}[theorem]{Remark}
\newtheorem{conjecture}[theorem]{Conjecture}

\numberwithin{equation}{section}

\newcommand{\Supp}{{\rm Supp}}

{\vskip-\lastskip\medskip
  \noindent
  {\em #1.}\enspace
  }%
{\qed\par\medskip
  }

\begin{document}
\title
[Hodge  Cohomology  Criteria  For Affine  Varieties]
{  Hodge  Cohomology  Criteria  For Affine  Varieties  }

\author[Jing  Zhang]{Jing  Zhang}

\begin{abstract}

We give several new criteria for a quasi-projective
variety to be affine. In particular, we prove that an algebraic manifold $Y$
  with dimension $n$  is affine if and only if  
$H^i(Y, \Omega^j_Y)=0$   for all $j\geq 0$,  $i>0$
and $\kappa(D, X)=n$, i.e., there are 
$n$ algebraically independent nonconstant regular functions on $Y$, 
  where $X$ is the smooth completion 
of $Y$, $D$ is the effective boundary divisor with support $X-Y$
and  $\Omega^j_Y$  is the sheaf of regular  $j$-forms on $Y$. 
This  proves  Mohan Kumar's affineness conjecture
for algebraic manifolds and     gives a  partial  answer
to   J.-P. Serre's Steinness question
\cite{36} in algebraic case since  the associated analytic  space of
an affine  variety  is  Stein
[15, Chapter VI, Proposition 3.1].

\end{abstract} 

\maketitle

\begin{center}
Department of Mathematics, University of 
Missouri, Columbia, MO
65211, USA
\end{center}
Email address: zhangj@math.missouri.edu

\begin{flushleft}
\end{flushleft}
\date{}
\section{Introduction}

 J.-P. Serre \cite{36} raised the following question in 1953: 
If   $Y$  is   a complex manifold
with   $H^i (Y, \Omega^j_Y)=0$   for all $j\geq 0$ and $i>0$, 
then what is $Y$? Is Y Stein?  Here  $\Omega^j_Y$
is  the sheaf  of  holomorphic 
$j$-forms  and  the cohomology is
${\check{\mbox{C}}}$ech  cohomology.  
A   complex   space  $Y$
 is Stein if and only if 
 $H^i(Y, G)=0$  for every analytic coherent sheaf $G$ on $Y$
 and all positive integers $i$. 
 For a holomorphic variety $Y$,
it is Stein if and only if  it is both holomorphically 
 convex and holomorphically separable  [13, Page 14]). We say 
that  $Y$ is  
holomorphically 
 convex if for any discrete sequence $\{y_n\}\subset Y$,
 there is a holomorphic function $f$ on $Y$ such that
 the supremum of the set $\{|f(y_n)|\}$ is 
 $\infty$. $Y$ is  
holomorphically  separable if for every pair $x,y\in Y$,
$x\neq y$,  there is a holomorphic function $f$ on $Y$ such that
$f(x)\neq f(y)$.

To investigate Serre's question,  
we are interested in the classification of  algebraic manifolds
(i.e., irreducible smooth algebraic  varieties  defined over $\Bbb{C}$)
with vanishing Hodge cohomology.  Since  the associated analytic variety of 
an  affine variety    is  Stein
[15, Chapter VI, Proposition  3.1],  it is natural to  ask
the algebraic  analogue  of  Serre's    question: 
if $Y$ is a smooth quasi-projective
variety with $H^i (Y, \Omega^j_Y)=0$   for all $j\geq 0$ and $i>0$
(from now  on   $\Omega^j_Y$
is  the sheaf  of  regular  
$j$-forms  and  the cohomology is still
${\check{\mbox{C}}}$ech  cohomology), 
then  is $Y$  affine? 
The answer is no even in surface case \cite{24}.
There are three types of surfaces and two of them 
are not affine since they have no nonconstant 
regular functions  [24, Lemma 1.8]. 
Curve case is trivial since $Y$ is not complete by Serre duality 
and any noncomplete curve is affine    
[15, Chapter II, Proposition 4.1]. 
We  have  classified  threefolds  with vanishing Hodge cohomology
in our previous papers [42, 43, 44, 45].
It again shows that $Y$ may not be affine.

 It is obvious that  
$H^i(Y, \Omega^j_Y)=0$   for all $j\geq 0$   and   $i>0$
is a necessary condition for  $Y$  being affine 
 since 
$\Omega^j_Y$   is coherent  [14, Chapter III, Section 5].
Conversely, since the vanishing Hodge cohomology is not 
a sufficient condition, a further question  is: 
What condition  should we add such that $Y$ is affine?  
This question is very interesting on its own because 
among algebraic varieties, affine varieties are basic, natural and important. 
 To study a nonaffine variety,  we can cover    it   by 
 affine varieties and examine the local data then glue  them 
 together to get global information by cohomology.

Now  we translate  our  question to a purely
algebraic  geometry question: 
When is a quasi-projective variety (or a scheme) affine? 
 There are many criteria  for  affineness. 
Here we only mention what we know. There may  be
other results we are not aware of.

  Serre's well-known criterion 
  says that a variety $Y$ is affine if and only if 
 $H^i(Y, F)=0$ for every coherent sheaf $F$ on $Y$
 and all positive integers $i$ [15, Chapter II, Section 1,
Theorem 1.1]. 
 Goodman and Hartshorne 
proved that $Y$ is affine if and only if 
$Y$ contains no complete curves and the dimension $h^1(Y, F)$
of the linear space 
$H^1(Y, F)$ is bounded for all coherent sheaf $F$ on $Y$ \cite{11}. 
Let $X$ be the completion of $Y$. 
Goodman also proved that $Y$ is affine if and only if 
after suitable blowing up
the closed subvariety on the  boundary  $X-Y$, the new boundary $X'-Y$
is the support of an ample divisor, where $X'\rightarrow X$
is the  blowing up with center in $X-Y$
[10; 15, Chapter 2, Theorem 6.1].  
Let $D$ be the effective boundary divisor with 
support $X-Y$,    then $Y$ is affine if $D$ is
ample. So  if we can show the ampleness of
$D$, $Y$ is affine. There are two important 
criteria for ampleness due to Nakai-Moishezon and 
Kleiman  
 [22; 23, Chapter 1, Section  1.5]. 
Another sufficient condition is that if 
$Y$ contains no complete curves and the linear system 
$|nD|$ is base point free, then $Y$ is  affine  [15, 
Chapter 2, Page 64].  Therefore we can apply base point free theorem
if we know the numerical condition of $D$ [ 32; 23, 
Chapter 3, Page 75, Theorem 3.3].
Neeman Proved that if $Y$
(the associated  scheme  of  $Y$)
 can be embedded in an affine 
scheme, then $Y$ is affine if and only if 
$ H^i(Y, {\mathcal{O}}_Y )=0 $  for all $i>0$  \cite{29}.

In higher dimension (at least, in our problem),
it is very hard to check  the ampleness of a big
(even big and nef)  divisor $D$
and  the base point freeness of the linear system 
$|nD|$. In \cite{45}, only  for a threefold 
and a prime divisor $D$, we show that  $|nD|$
is base point free for all sufficiently large $n$.

The theory of $D$-dimension and Iitaka  fibration 
is widely used  in the classification   of algebraic  varieties
[17-20, 26, 39].  
The notion of $D$-dimension is due to Iitaka \cite{17}:
If for all $m>0$, $H^0(X, {\mathcal{O}}_X(mD))=0$, then 
the $D$-dimension
$\kappa(D, X)=-\infty$. Otherwise, 
$$\kappa(D, X)=tr.deg_{\Bbb{C}}\oplus_{m\geq 0}
H^0(X, {\mathcal{O}}_X(mD))-1.$$
 
Based on the classification of surfaces \cite{24}, 
 for threefolds, Mohan Kumar conjectured 
that $Y$ is affine if and only if  
$H^i(Y, \Omega^j_Y)=0$   for all $j\geq 0$,  $i>0$
and $\kappa(D, X)=3$, i.e., there are 3
algebraically independent nonconstant regular functions on $Y$,
where $X$ is the smooth completion 
of $Y$ and $D$ is the boundary divisor with support $X-Y$. 
If $Y$ is a smooth threefold  and  $X-Y$  is a smooth projective
surface, 
we have proved this conjecture  
in \cite{45}.   In fact, we can  drop the 
smoothness assumption  when $Y$ is a threefold.

We work over  the complex  number field  $\Bbb{C}$.

\begin{theorem}If $Y$ is a  quasi-projective 
threefold, then $Y$ is affine if and only if 
 $H^i(Y, \Omega^j_Y)=0$   for all $j\geq 0$, $i>0$
and $\kappa(D, X)=3$, where $X$ is a completion of $Y$
and $D$ is an effective boundary divisor with support $X-Y$. 
\end{theorem}

In higher dimension, we have the following conjecture.

\begin{conjecture} An algebraic manifold $Y$ of dimension $n$
 is affine if and only if   
$H^i(Y, \Omega^j_Y)=0$   for all $j\geq 0$,  $i>0$
and $\kappa(D, X)=n$. 
\end{conjecture}

In fact, we can  prove  a 
  more general theorem 
since  we do not require the smoothness of the variety.

\begin{definition}
 We say that $Y$ satisfies the condition A if 
the following three conditions hold

(1) the boundary $X-Y$
is connected;

(2)  $X-Y$ is of  pure
codimension 1 in $X$;

(3) every  closed subvariety $Z$  of  codimension 1 in  $Y$   is affine.

\end{definition}

 \begin{theorem} If $Y$ is an  irreducible   
 quasi-projective variety 
     with dimension $n\geq 1$,
then $Y$ is affine if and only if the following three conditions 
are satisfied:

(1) $Y$ satisfied condition A;

(2) $\kappa(D, X)\geq 1$;

(3) $ H^i(Y, {\mathcal{O}}_Y )=0 $  for all $i>0$. 

\end{theorem} 

\begin{theorem} If $Y$ is a smooth  quasi-projective
variety with dimension $n$, then $Y$ is affine if and only if  
$H^i(Y, \Omega^j_Y)=0$   for all $j\geq 0$,  $i>0$
and $\kappa(D, X)=n$.
\end{theorem}

By the   result  in  \cite{42}, if $Y$ is an algebraic manifold with 
$H^i(Y, \Omega^j_Y)=0$   for all $j\geq 0$,  $i>0$
and $\kappa(D, X)\geq 1$, then there is a surjective morphism from 
$Y$ to a smooth affine curve  $C$ such that a general fibre is smooth 
and irreducible.   In   \cite{42}, we asked 
whether  $Y$ is affine  if  
 the general fibre
is affine.  By the above theorem, it is easy to 
see that the answer is yes.

\begin{theorem} If $Y$ is a smooth  quasi-projective
variety with dimension $n$, then $Y$ is affine if and only if  
$H^i(Y, \Omega^j_Y)=0$   for all $j\geq 0$,  $i>0$, 
 $\kappa(D, X)\geq 1$ and the general fibre of $Y\rightarrow C$ is affine.
\end{theorem}

\begin{theorem} If $Y$ is a smooth  quasi-projective
variety with dimension $n$, then $Y$ is affine if and only if  
$H^i(Y, \Omega^j_Y)=0$   for all $j\geq 0$,  $i>0$, 
 $\kappa(D, X)\geq 1$ and   every   closed  subvariety 
of codimension 1  is affine.
\end{theorem}

The following corollary is a partial answer to J.-P. Serre's 
Steinness question. 
\begin{corollary} If an n-dimensional  algebraic manifold  $Y$
with  $H^i(Y, \Omega^j_Y)=0$   for all $j\geq 0$ and  $i>0$
has  $D$-dimension $n$, then $Y$ is   Stein. 
\end{corollary} 

\begin{definition}An algebraic variety $Y$ is regularly separable if 
for any two distinct points $y_1$ and $y_2$ on 
$Y$, there is a regular
function $f$ on $Y$ such that $f(y_1)\neq f(y_2)$.
\end{definition}

\begin{corollary} If an   algebraic manifold  $Y$ of 
dimension $n$  satisfies 
 $H^i(Y, \Omega^j_Y)=0$   for all $j\geq 0$ and  $i>0$, then 
the following  conditions are equivalent

(1) $Y$ is affine;

(2) $Y$  is regularly separable;

(3) any closed   subvariety of codimension 1
in $Y$ is affine; 

(4) $\kappa(D, X)=n$.

\end{corollary}

This paper is organized as follows. 
We will   give two preliminary  lemmas and review  some  known
  results   in Section 2 
and 
 prove the above theorems in  Section 3. 
The main idea of  our proof  is to  show  that $Y$ is regularly separable,
i.e.,  for any two distinct points $y_1$ and $y_2$, there is
a regular function $H$ on $Y$ such that $H(y_1)\neq H(y_2)$. 
Then we have an injective  birational  
 morphism from the associated scheme of $Y$ 
to 
Spec$\Gamma(Y, {\mathcal{O}}_Y)$.   Therefore we can  apply
Zariski's Main  Theorem  [27, Chapter III, Section 9]
  and    Neeman's 
result  [29].\\

\noindent 
{\bf{Acknowledgments}}  \quad I would like to thank 
the following 
      professors for  helpful discussions: Steven Dale Cutkosky,
Dan  Edidin, 
N.Mohan Kumar, Zhenbo Qin  and Qi Zhang. In particular, I thank 
Mohan Kumar who suggested me to construct an embedding 
then use Neeman's result.

\section{Preliminary  Lemmas}

The proof of Lemma 2.1 is essentially  due to  Mohan Kumar
[24]. We include a proof here for completeness.

\begin{lemma} Let $Y$ be a  smooth   quasi-projective
   variety  with 
 $H^i(Y,\Omega^j_Y)=0$
for every $j\geq 0$ and $i>0$. Let 
$X$ be the projective variety containing  
 $Y$, then
$X-Y$  is connected. 
\end{lemma}

$Proof$.
If $X-Y=Z$ is not connected, write $Z=Z_1+Z_2$, 
where both $Z_1$ 
and $Z_2$ are closed in $Z$, nonempty and 
 $Z_1\cap    
     Z_2=\emptyset.$ 
Let dim$Y$=n. We have a long exact sequence of local cohomology
[14, Page 212]
$$0=H^{n-1}(Y,\Omega^n_Y)\longrightarrow 
      H^n_Z(X,\Omega^n_X)\longrightarrow 
 H^n(X,\Omega^n_X)\longrightarrow  H^n(Y,\Omega^n_Y)
=0.$$
By Serre duality
$$ H^n_Z(X,\Omega^n_X)=H^n(X,\Omega^n_X)=H^0(X, {\mathcal{O}}_X)=\Bbb{C}.$$ 
But by Mayer-Vietoris sequence   [14, Page 212],
$$ H^n_Z(X,\Omega^n_X)\cong  H^n_{Z_1}(X,\Omega^n_X)\oplus 
 H^n_{Z_2}(X,\Omega^n_X).$$ 
Both summands are at least one dimensional since 
$$ H^n_{Z_i}(X,\Omega^n_X)\longrightarrow  H^n(X,\Omega^n_X)={\Bbb{C}}
 \longrightarrow    H^n(X-Z_i,\Omega^n_X)=0,$$
where the last cohomology vanishes   because  $X-Z_i$ is not complete
[21, 37]. 
This is a contradiction.   
\begin{flushright}
 Q.E.D. 
\end{flushright}

\begin{remark}  The above proof also works for 
    singular varieties because for singular complete varieties,
we have Grothendieck  duality 
[1, Chapter I, 1.3]: $H^n(X, \omega_X)\cong  \Bbb{C}$.    
\end{remark}

Recall an  equivalent definition  of  
$D$-dimension.  If  $h^0(X, {\mathcal{O}}_X(mD))>  0$
for  some  $m\in {\Bbb{Z}}$  and $X$  is normal, 
choose a basis $\{f_0, f_1, \cdot \cdot\cdot, f_n\}$
     of the linear space 
     $H^0(X, {\mathcal{O}}_X(mD))$, it defines a rational 
     map 
     $\Phi _{|mD|}$
     from $X$ to the projective space 
     ${\Bbb{P}}^n$ by sending a point $x$ on $X$ to
     $(f_0(x), f_1(x), \cdot \cdot\cdot, f_n(x))$ in ${\Bbb{P}}^n$. 
By definition of $D$-dimension [39, Definition 5.1],
$$ \kappa (D, X)= \max_m\{\dim (\Phi _{|mD|}(X))\}. 
      $$

\begin{lemma} Let $X$, $Y$ be as above and $n$ be the dimension of $Y$.
Then 
 $Y$ contains no complete curves if $\kappa(D, X)=n$, where 
$D$ is an effective divisor with support $X-Y$. 
 \end{lemma}
       $Proof$. We will prove the claim by induction
on the dimension of $Y$. The claim is trivial if $Y$ 
is a curve since $Y$ is not complete  by Serre 
duality. The surface case was proved in \cite{24}.  Assume
that the claim holds for  every $(n-1)$-dimensional 
 variety with the same property.
Suppose now  dim$Y=n$. If $Y$ has a connected complete curve $C$, 
pick a point $p$ on $C$. Let $Z$  be a smooth prime principal divisor  
passing through $p$ (for the existence of $Z$, see Step 2 of proof 
of Theorem 1.5 in next section). 
Let $f\in H^0(Y,{\mathcal{O}}_Y)$ be the defining function of 
$Z$, then we have the following short exact sequence 
$$ 0\longrightarrow 
 \Omega^j_Y
\longrightarrow 
 \Omega^j_Y
\longrightarrow 
 {\Omega^j_Y}|_Z
\longrightarrow 
0,
$$
where the first map is  defined by $f$. 
Thus for all $i>0$ and $j\geq 0$,  we have 
$$ H^i(Z,  {\Omega^j_Y}|_Z)=0.
$$
In particulsr,  $H^i(Z, {\mathcal{O}}_Z)=0$. 
From the exact sequence [14, Chapter II, Theorem 8.17; 
12,   Page  157],
$$0\longrightarrow 
\Omega^{j-1}_Z\longrightarrow
  \Omega^j_Y|_Z
\longrightarrow 
\Omega^j_Z
\longrightarrow 
0, 
$$
we have 
$$ H^i(Z,  \Omega^j_Z)=0
$$
for all $i>0$ and $j\geq 0$. 
By the inductive assumption, $Z$ contains no complete curves.

Consider the morphism from $Y$ to $\Bbb{C}$
$$ f: Y\longrightarrow  \Bbb{C}.
$$
We know $p\in Z = f^{-1}(0)$, a prime   principal  divisor 
defined  by   $f$.  Since  $C$
is connected,   $C\subset   Z  $.  This is in
contradiction with the fact that 
$Z$   has no complete curves. 

The proof is completed.

\begin{flushright}
 Q.E.D. 
\end{flushright}

\begin{theorem}{\bf[Fujita]} Let  $D$ be an effective  $\Bbb{Q}$-divisor
on a normal projective  surface $X$. Then  there exists a unique 
decomposition 
$$D=P+N$$
satisfying the following conditions:

(1) N is an effective $\Bbb{Q}$-divisor and either N=0 or the intersection 
matrix of the irreducible components of N is negative definite;

(2) P is a nef $\Bbb{Q}$-divisor and the intersection of P 
with each irreducible 
component of N is zero. 
\end{theorem}

\begin{lemma}{\bf[Sakai]}   Let $D=P+N$ be the Zariski decomposition 
of an effective divisor $D$ on a normal complete surface  $X$, 
then  $\kappa(D, X)=2$ if and only if $P^2>0$.
\end{lemma}

\begin{theorem}{\bf[Iitaka]}  Let $X$ be a normal projective variety and let 
$D$ be an effective divisor on $X$. Then  there exist two positive
numbers $\alpha$ and $\beta$ such that for all sufficiently 
large $n$ we have
$$  \alpha n^{\kappa(D, X)}
\leq h^0(X, {\mathcal{O}}_X(nD))
\leq \beta n^{\kappa(D, X)}. 
$$
\end{theorem}
 We say that $D$
is a big  divisor if   $h^0(X, {\mathcal{O}}_X(nD))\geq \alpha n^d$, where 
$d$ is the dimension of $X$.

A fibre space  is a  proper surjective   morphism
$f: V\rightarrow  W$  between  two varieties  $V$ and $W$
such that  the  general   fibre  is  connected.  We know the following two 
ways to calculate  the $D$-dimension of a  variety. 

\begin{theorem}{\bf[Iitaka]}  Let  $f: V\rightarrow  W$
be a  surjective  morphism  of 
two varieties  and let  $D$ be  a Cartier divisor  on $W$,
then we have 
$$  \kappa(f^*D,   V)=
\kappa(D, W). 
$$
\end{theorem}

\begin{theorem}{\bf[Iitaka]}  Let  $f: V\rightarrow  W$
be a   fibre space  from  a complete non-singular variety $V$
to  a  variety $W$.  If  $D$  is a divisor on  $V$, then 
there is  an open dense subset (in complex  topology) $U$
in $W$  such that  for every  point $u\in  U$
$$  \kappa(D,   V)\leq  \kappa(D_u,   V_u) + dim W,
$$
where  $V_u=f^{-1}(u)$ and  $D_u=D|_{V_u}$.
\end{theorem}   
For the proof of  Iitaka's   theorems, see  [18, Lecture 3]
or  [39,  Chapter II, Section 5].

We will frequently change the coefficients of the effective 
divisor $D$  or blow up a closed subvariety in $D$  without changing
the $D$-dimension because of the following 
 two  properties [17; 39,  Chapter II, Section 5].  

Let  $f$: $X'\rightarrow X$
        be a surjective morphism between  two complete
        varieties $X'$ and $X$, let $D$ be a divisor on 
        $X$ and $E$ an effective divisor on $X'$ such that
        codim$f(E)\geq 2$, then  
        $$\kappa (f^{-1}(D)+E, X')=\kappa (D, X), 
        $$
where  $f^{-1}(D)$   is the reduced  transform of $D$,
defined to be
$f^{-1}(D)=\sum D_i$, $D_i$'s are the   irreducible components
of  $D$. 
       The second property of $D$-dimension is that 
       it does not depend on the coefficients of $D$ 
       under  a mild condition which is    true in our case
       since we always choose effective boundary divisor $D$
       with normal crossings.
       Let $D_1$, $D_2$, $\cdot$$\cdot$$\cdot$, $D_n$
       be  any divisor on $X$ such that for every $i$, $0\leq i \leq n$, 
       $\kappa (D_i, X)\geq 0$, then for   integers 
       $p_1> 0,\cdot\cdot\cdot$, $p_n>0$,  we have   [20, Section 5]
       $$\kappa(D_1+\cdot\cdot\cdot+D_n,X)=
       \kappa(p_1D_1+\cdot\cdot\cdot+p_nD_n,X).
       $$
      In particular, if $D_i$'s are irreducible components of
      $D$ and $D$ is effective, then we can change its coefficients
      and do not change the $D$-dimension.

For a projective manifold $M$, let $L$
be a line bundle on $M$, then  there is a  Cartier divisor 
$D$ determined by $L$. We define $\kappa(L,M)=\kappa(D,M)$.

\begin{lemma}[{\bf Fujita}] Let 
$M$ and  $S$  be two projective manifolds.  
 Let $\pi: M\rightarrow S$
  be a  fibre space and let $L$ and $H$ be  line bundles  on $M$
  and $S$ respectively. Suppose that  
  $\kappa (H, S)=\dim S$ and that
  $\kappa (aL-b\pi^*(H))\geq 0$
for certain  positive integers $a$, $b$.  Then
$\kappa(L,M)=\kappa (L|_F, F)+ \kappa (H, S)$
for a general  fibre $F$  of  $\pi$.   
\end{lemma}

Let ${\Bbb{C}}(X)$ be the function field of $X$. Let 
$$  R(X, D)=\oplus_{\gamma\geq 0}H^0(X, {\mathcal{O}}_X(\gamma D)) 
$$
be the graded $\Bbb{C}$-domain  and $R^*\subset R$
the  multiplicative  subset  of all nonzero
homogeneous  elements. Then the quotient ring $R^{*-1}R$  
is a graded  $\Bbb{C}$-domain  and its degree 0 part 
$(R^{*-1}R)_0$  is a  field  denoted  by  $Q((R))$, or $Q((X, D))$, i.e., 
$$ Q((R))=Q((X, D))=(R^{*-1}R)_0.
$$  
The proof of the following   lemma can be found in \cite{26}. 
\begin{lemma} Let $X$ be  normal proper  over an algebraically closed field
$k$. If there is an $m_0>0$  such that  for all $m>m_0$, 
$h^0(X, {\mathcal{O}}_X(mD))>0$, then 
$$ {\Bbb{C}}(\Phi _{|mD|}(X)) = Q((X, D)).
$$  
In particular, if  $\kappa(D, X)=$\mbox{dim}X, then  
$\Phi _{|mD|}$ is birational for all $m\gg 0$.
\end{lemma}

There are many versions of Zariski's Main Theorems.
The proof of the following Zariski's Main Theorem can 
be found in   [27, Chapter III, Section 9].

\begin{theorem} {\bf[Zariski's Main Theorem]} 
Let $X$ be a normal variety  over a field $k$ and
let $f:X'\longrightarrow X$
be a birational morphism with finite 
fibre from a variety $X'$ to $X$.
Then   $f$  is an isomorphism of  $X'$
with an open subset $U\subset X$. 
\end{theorem}

\begin{theorem}{\bf [Neeman]} Let  $X=$${\mbox{Spec}}$A be a scheme,
$U\subset X$  a quasi-compact  Zariski  open  subset. Then
$U$  is affine  if and only if  $H^i(U, {\mathcal{O}}_U)=0$
for $i\geq 1$.
\end{theorem}

\section{Proof of the Theorems}

Recall our notation:  $Y$ is an open subset of 
a projective variety $X$ and $D$ is the effective boundary divisor
 with support $X-Y$.  We may assume that the boundary 
divisor $D$ has simple normal crossings by further blowing up suitable
closed subvariety  of $X-Y$.

{\bf Proof of Theorem 1.1.} One direction is trivial. 
If $Y$ is affine, then
 $H^i(Y, \Omega^j_Y)=0$   for all $j\geq 0$,  $i>0$
since $\Omega^j_Y$ is a coherent  sheaf  [14, Chapter 3, Theorem 3.7];
 the affineness of $Y$ also implies that  $D$ is big since after
 further blowing up the boundary
$X-Y$, it is the support of an ample divisor 
[ 15, Chapter 2, Section 6, Theorem 6.1].
Therefore  $\kappa(D, X)=3$.

Conversely, 
we need to prove  that   $Y$ is affine
if $H^i(Y, \Omega^j_Y)=0$   for all $j\geq 0$,  $i>0$
and $\kappa(D, X)=3$.
It is sufficient to prove the theorem 
for  normal  varieties. The reason is the following.

 By Chevally's theorem
[14, Page 222;  15, Chapter 2, Page 63],
$Y$ is affine if and only if
its normalization $Y'$ is affine since we have a finite morphism from $Y'$
to $Y$.  Because a finite morphism is an affine morphism  [14, Page 128],
 $Y'$ also satisfies the vanishing cohomology 
$H^i(Y, \Omega^j_Y)=0$   for all $j\geq 0$ and $i>0$
 [14, page 
222].   By the  definition and property  of $D$-dimension,
the new tripple  $(Y', X', D')$ after the normalization also satisfies 
the hypothesis of  Theorem 1.1  [17].  So we may assume that $Y$ is normal.

We divide our proof into several steps.

{\bf Step 1}.  Every principal divisor $Z=\{f=0\}$ defined by a 
nonconstant regular 
function $f\in H^0(Y, {\mathcal{O}}_Y)$ satisfies  
$H^i(Z, \Omega^j_Z)=0$   for all $j\geq 0$,  $i>0$.

$Proof$.
For $j=0,1,2,$ the function $f$ gives a short exact sequence
$$ 0\longrightarrow 
\Omega^j_Y
\longrightarrow 
\Omega^j_Y
\longrightarrow 
\Omega^j_Y|_Z
\longrightarrow 
0, 
$$
where the first map is defined by $f.$
Since  $H^i(Y, \Omega^j_Y)=0$ 
 for every $i>0$ and $j\geq 0$,
$H^i(Z, \Omega^j_Y|_Z)=0$. In particular, 
$H^i(Z, {\mathcal{O}}_Z)=0$ for all $i>0$. 
We have an exact sequence   [14, Chapter 2, Section 8],
$$0\longrightarrow 
A
\longrightarrow 
\Omega^j_Y|_Z
\longrightarrow 
\Omega^j_Z
\longrightarrow 
0, 
$$ 
where $A$ is a coherent sheaf on $Z$. Since $Z$ is not complete,
 $H^2(Z, A)=0$   [21]. Thus we  have
$H^i(Z, \Omega^j_Z)=0$   for all $j\geq 0$,  $i>0$. \\

\begin{remark} If the dimension of $Y$ is higher than 3,
in the last exact sequence, in  general  
 $H^2(Z, A)$ does not vanish. Only when $Z$ is smooth, for all 
$i>0$, 
 $H^i(Z, A)=0$ if  dim$Y>$3.  
\end{remark}

{\bf Step 2}. Let $\bar{Z}=X_0\subset X$ be a projective surface 
in $X$ such that $Z$ is an open subset of $\bar{Z}$.
Let $D_0$ be the boundary divisor with support $\bar{Z}-Z$.
 If $Z$ is a surface with 
$H^i(Z, \Omega^j_Z)=0$   for all $j\geq 0$,  $i>0$
and $\kappa(D_0, X_0)=2$, then $Z$ is affine.

$Proof.$ We proved this result in   [44].
We include a proof here for completeness. 

By the same proof of Lemma 1.1 and Lemma 1.4 in   [24],
 we know that $Z$ contains no complete curves
and the boundary $D_0$ is  connected.  We will prove that 
$Z$ is affine   if   it satisfies the following three conditions

(1) $Z$ contains no complete curves;

(2) $\bar{Z}-Z$ is connected;

(3) $\kappa(D_0, X_0)=2$.

The idea of proof  is to   show  that 
the boundary $\bar{Z}-Z$ is  the support of an ample divisor  $P$
on $\bar{Z}$. 

Again we may assume that $Z$ is normal. 
For a normal surface $X$, the intersection theory  is due to
Mumford \cite{28}. For any effective divisor on a complete normal surface,
we have Zariski decomposition \cite{34}.

Write the Zariski decomposition  $D=P+N$, where 
$N$ is negative  definite, $P$  is effective and nef and 
any prime component of  
$N$  does not intersect $P$   \cite{41, 34}. We may assume that 
both $P$ and $N$ are integral by multiplying a positive 
integer to  $D$. 
Let  $\Supp{D}=\{D_1, D_2, \cdot\cdot\cdot, D_n \}=X-Y$.
Since   $\kappa (D, X)=2$,  $P^2>0$   [3 Corollary 14.18]. 
First we claim that $\Supp{P}=\Supp{D}=X-Y$.

If 
$\Supp{P}\neq X-Y$, then there is a prime component   $D_1$
in $X-Y$ such that  $P\cdot D_1>0$ and  $D_1$ is not a component of 
$P$  since $X-Y$ is connected.   Let 
$$Q=mP+D_1,$$
 where $m$ is  a  big  positive integer.
Then  $Q$ is an effective divisor and $\Supp{Q}=\Supp{P}\cup D_1$. 
We may choose $m$ such  that  
$$ Q^2=m^2P^2+2mP\cdot  D_1+ D_1^2>0.  
$$
For every prime  component  $E$ in $P$, since $P$ is nef and 
$D_1$ is not contained in $\Supp {P}$, we  can choose
sufficiently large $m$ such  that
$$ Q\cdot E=mP\cdot  E+D_1\cdot E\geq 0, \quad \quad 
D_1\cdot  Q=mD_1\cdot P +D_1^2 >0.
$$
Thus we get a new effective divisor $Q$ 
such  that  $Q$ is nef and $Q^2>0$. We may replace $P$ by $Q$
and still call it $P$. 
By finitely many such   replacements, we 
can get an effective nef  divisor $P$     such that 
$P^2>0$ and $\Supp{P}=\Supp{D}=X-Y$.

We  claim  that the boundary $X-Y$ is the support of an ample divisor. 
In fact, the following three conditions imply  the   ampleness: 
(1) $X-Y$ is connected; (2) $Y$ contains no complete curves;  
(3) There is an effective nef divisor $P$ with supp$P=X-Y$ and $P^2>0$.

If $P$ is not ample, then there is an irreducible 
curve $C$ in $X$ such that
$P\cdot C=0$ by Nakai-Moizshon's ampleness criterion
[14, Chapter V, Theorem 1.10]. Since $Y$ has no complete curves, 
$C$ must be one of the $D_i's$. Rearrange the order, we may assume 
$D_i\cdot P=0$ for $i=1,2,..., r$ and $D_j\cdot P>0$
for $j= r+1,..., n$.  
Write 
$$P=\sum_{i=1}^ra_iD_i+\sum_{j=r+1}^nb_jD_j=A+B,$$
where
$A=\sum_{i=1}^ra_iD_i$, $B=\sum_{j=r+1}^nb_jD_j=A+B$. 
Then for $i=1,...,r$,
$$0=P\cdot D_i = A\cdot  D_i+B\cdot D_i.$$
Since $D_i$ is not a component of $B$ for $i=1,...,r$,
$B\cdot D_i\geq 0$. So $A\cdot  D_i\leq  0$  for every 
$i=1,..., r$. Thus the intersection matrix 
$[D_s\cdot D_t]_{1\leq s, t\leq r}$
is  negative  semi-definite \cite{2}.
Since $A\cup B=X-Y$ is connected, there is at least one 
component of $A$, say, $D_{i_o}$, such that 
$D_{i_o}\cdot B>0$. Hence 
 $D_{i_o}\cdot A<0$. This implies that 
the intersection matrix 
$[D_s\cdot D_t]_{1\leq s, t\leq r}$
is  negative  definite \cite{2}.
Therefore there is an effective divisor 
$E=\sum_{i=1}^r\alpha _iD_i$  such that 
$E\cdot D_i<0 $ for all 
$i=1,...,r$  \cite{2}.

So there are positive numbers 
$\alpha_i$, $i=1,..., r$, such that for every $i$, $E\cdot D_i< 0$,
where $E= \sum_{i=1}^r\alpha_iD_i$.  Let $P_1=mP-E$, $m\gg 0$, then
$P_1^2>0$, $P_1$ is nef  and if $1\leq  i\leq r$, 
$$P_1\cdot  D_i=-E\cdot D_i>0.$$
If $r+1 \leq j\leq n$, then  choose sufficiently large 
$m$ such that
$$  P_1\cdot  D_j=mP\cdot D_j-E\cdot D_j>0.
$$
Thus $P_1$ is an effective ample divisor with support $X-Y$. Replace 
$P$ by $P_1$, we    have shown  that $X-Y$ is the support of an ample divisor 
$P$.  
Therefore $Y$ is an affine surface. \\

{\bf Step 3}.  An irreducible principal divisor 
$Z=\{f=0, f\in H^0(Y, {\mathcal{O}}_Y)\}$
is affine. 

$Proof$. In Step 1, we showed  that 
$H^i(Z, \Omega^j_Z)=0$   for all $j\geq 0$,  $i>0$.
We will prove $\kappa(D_0, X_0)=2$ then apply Step 2,
where  $\bar{Z}$ 
is the closure  of  $Z$
 is $X$,  $X_0=\bar{Z}\subset X$, supp$D=\bar{Z}-Z$.

 By \cite{42}, the regular function $f$  determines the
 following commutative diagram
\[
  \begin{array}{ccc}
    Y                           &
     {\hookrightarrow} &
    X                                 \\
    \Big\downarrow\vcenter{%
        \rlap{$\scriptstyle{f|_Y}$}}              &  &
    \Big\downarrow\vcenter{%
       \rlap{$\scriptstyle{f}$}}      \\
C        & \hookrightarrow &
\bar{C}
\end{array}
\]
where $C$ is a smooth  affine curve 
embedded in a smooth projective 
curve $\bar{C}$, and $f$ is proper and surjective, 
every fibre of $f$ over $\bar{C}$
 is connected, 
 general fibre  is  smooth. Also general fibre 
 of $f|_Y$ is connected  and smooth. In particular,  $Z=\{f=0\}$
is one fibre in $Y$. Let $\bar{Z}=f^{-1}(0)$ be the inverse image of $0$
in $X$, then $Z=\bar{Z}\cap Y$.

Let  $\pi :X'\rightarrow X$ be the proper surjective morphism from
a smooth projective variety $X'$ to $X$ with connected fibres
($\pi$ is birational and generically finite), 
then we have a new fibre
space  $f': X'\rightarrow \bar{C}$ with the same property as 
the original fibre  space $f: X\rightarrow   \bar{C}$. 
By  [39, Theorem 5.1], 
$$ \kappa(\pi ^*D, X')=\kappa(D, X)=3.
$$
By  [39, Theorem 5.11], 
for a general fibre $X'_t=f'^{-1}(t)$ in $X'$, $t\in  C$, 
we have 
$$  3\leq  \kappa(f'^*D, X') \leq 
\kappa(f'^*D|_{X'_t}, X'_t)+1, 
$$
where $X'_t=f'^{-1}(t)$ is the the general fibre in $X'$.
So $\kappa(f'^*D|_{X'_t}, X'_t)=2$. 
Apply    [39,  Theorem 5.13]   to the fibre $X'_t$, we have
$$  \kappa(D_t, X_t)= \kappa(f'^*D|_{X'_t}, X'_t)= 2.
$$ 
By upper semi-continuity theorem [14,
Chapter 3, Section 12]  or [39, Chapter 1], 
$\kappa(D_0, X_0)=2$. By Step 2, $Z$ is affine. \\

{\bf Step 4}.  
For every point $y\in Y$, there is a principal  
divisor $Z$ passing through $y$ and the closure $\bar{Z}$ of $Z$ in $X$
is  connected.

$Proof.$  The following construction of fibre space 
is due to Ueno  [40, Page 46].

Since $\kappa(D, X)=3$, 
let  $\Phi _{|nD|}$ be the rational map from $X$ to  ${\Bbb{P}}^N$ 
defined  by $\{\phi_0, \phi_1, \cdot\cdot\cdot, \phi_N\}$, a basis of 
$H^0(X, {\mathcal{O}}_X(nD))$, where  $N$=dim$|nD|$ and dim(im($\Phi _{|nD|}))
\geq 2$.  Choose three hyperplane divisors $H_a$, $H_b$, and $H_c$ 
such that the rational functions $\eta_1$ and $\eta_2$ on $X$ 
induced by rational functions 
$$\frac{a_0X_0+a_1X_1+\cdot\cdot\cdot+a_NX_N}
{b_0X_0+b_1X_1+\cdot\cdot\cdot+b_NX_N}
$$
and 
$$\frac{c_0X_0+c_1X_1+\cdot\cdot\cdot+c_NX_N}
{b_0X_0+b_1X_1+\cdot\cdot\cdot+b_NX_N}
$$
on    ${\Bbb{P}}^N$ 
are algebraically independent, where $a=(a_0, a_1, \cdot\cdot\cdot, a_N)$, 
$b=(b_0, b_1, \cdot\cdot\cdot, b_N)$,
$c=(c_0, c_1, \cdot\cdot\cdot, c_N)$,
$a, b, c\in {\Bbb{C}}^{N+1}.$ By the rational map 
$\Phi _{|nD|}$, we can consider that $\eta_1$ and $\eta_2$ 
are elements of ${\Bbb{C}}(X)$. By Zariski's lemma   [16, Chapter X, 
Section 13, Theorem 1, Page 78], there exists a constant $d$ such that 
the field 
${\Bbb{C}}(\eta_1 +d\eta_2)$ is algebraically closed in ${\Bbb{C}}(X)$. 
Define a rational map $f$ from $X$ to  ${\Bbb{P}}^1$ by sending points $x$ 
in $X$ to 
$(1, \eta_1(x) +d\eta_2(x) )$ in  ${\Bbb{P}}^1.$  
We can choose 
 $\eta_1$ and $\eta_2$, 
such that  
$\eta_1 +d\eta_2$ 
only has poles in $D$, that is, when restricted to $Y$, $f$ is a morphism. 
By  [39, Corollary 1.10],  we have diagram 
\[
  \begin{array}{ccc}
    Y                           &
     {\hookrightarrow} &
    X                                 \\
    \Big\downarrow\vcenter{%
        \rlap{$\scriptstyle{f|_Y}$}}              &  &
    \Big\downarrow\vcenter{%
       \rlap{$\scriptstyle{f}$}}      \\
C        & \hookrightarrow &
{\Bbb{P}}^1,
\end{array}
\]
where  $f$ is proper and surjective, 
every fibre of $f$ over $C$
 is connected, 
 general fibre  is  smooth. Also general fibre 
 of $f|_Y$ is connected  and smooth.

For a point $y\in Y$, let $f(y)=a$, then 
$y\in S_a=f^{-1}(a)\cap Y$. By \cite{42}, 
$S_a$ satisfies $H^i(S_a, \Omega^j_{S_a})=0$
for all $i>0$ and $j\geq 0$. So $S_a$ contains no complete 
curves. 
Let $X_a=f^{-1}(a)$, then  the boundary $X_a-S_a$ is connected.
By the same argument as in the Step 3, we claim that
$\kappa(D_a, X_a)=2$, where $D_a$ is the effective 
boundary divisor supported in  $X_a-S_a$. 
By Step 2, $S_a$ is affine. Let $Z=S_a$, we are done.

{\bf Step 5}.
There is an injective morphism from  the associate scheme of  
$Y$
 to  {\mbox{Spec}}$\Gamma(Y, {\mathcal{O}}_Y)$.  

$Proof.$  First let $P$ be a point on the quasi-projective variety
 $Y$, then it is  a closed point on 
the associated scheme (we still call it $Y$) of $Y$. 
Let ${{\mathcal{M}}_P}$ be the maximal ideal of 
the local ring ${\mathcal{O}}_P$. Let 
$n=$dim$_{\Bbb{C}}{{\mathcal{M}}_P}/{{\mathcal{M}}_P}^2$,
then there are $n$ regular functions $f_1$,..., $f_n$
on an open subset containing 
$P$ such that
$${{\mathcal{M}}_P}=(f_1, ..., f_n).$$ 
Since $D$ is big, we claim that 
 each $f_i$   can be   extended  to a regular function on $Y$
(See Step 1, proof of Theorem 1.5).
 Then 
these $n$ regular functions give a 
maximal ideal of $\Gamma(Y, {\mathcal{O}}_Y).$

For an irreducible  closed subset 
$A$  of $Y$, let  $f_1$,..., $f_r$ be the local defining functions of
$A$, then we may assume that all of them are regular on $Y$.
So the ideal $I=(f_1,..., f_n)$  is a  prime ideal
of $\Gamma(Y, {\mathcal{O}}_Y)$. 
Naturally we can define  a map from $Y$ to 
Spec$\Gamma(Y, {\mathcal{O}}_Y)$ by sending a closed point to 
the corresponding maximal ideal in $\Gamma(Y, {\mathcal{O}}_Y)$ 
and a nonclosed  point  on the scheme $Y$ to the associated 
prime ideal in  $\Gamma(Y, {\mathcal{O}}_Y)$. 
Let $g$ be the morphism defined    by   this  map.  We will prove its 
injectivity.

  $g$ is injective if and only if 
for any two distinct points $y_1$ and $y_2$ in $Y$, 
there is a regular function $f$ on $Y$
separating these two points, i.e., $f(y_1)\neq f(y_2)$.
To see this, let $A$ and $B$ be two distinct irreducible closed   
subsets of $Y$ (as a quasi-projective variety).
There is a point $p$  in one of them, say $A$, such that $p$
in not a point of $B$. Then there is a regular function $f$
on $Y$ such that $f(y_1)=0$ but $f(y_2)=1$. Therefore $f$ is one of the 
defining functions of $A$ but not $B$. Write the 
prime ideals defining  $A$   and $B$   as follows
$$  g(A)=(f, f_1, ... f_m), \quad  g(B)=(g_1, ..., g_n).
$$
Then $f\not\in  g(B)$.  So
the images $g(A)\neq g(B)$.
This   shows the  injectivity of  $g$.

Let $y_1$ and $y_2$ be two distinct points in $Y$ and  
consider their images in the fibre space constructed in Step 4
\[
  \begin{array}{ccc}
    Y                           &
     {\hookrightarrow} &
    X                                 \\
    \Big\downarrow\vcenter{%
        \rlap{$\scriptstyle{f|_Y}$}}              &  &
    \Big\downarrow\vcenter{%
       \rlap{$\scriptstyle{f}$}}      \\
C        & \hookrightarrow &
{\Bbb{P}}^1.
\end{array}
\]
If $f(y_1)\neq f(y_2)$, then $g(y_1)\neq g(y_2)$. 
Assume  $f(y_1) = f(y_2)=a\in C=f(Y)\subset \Bbb{C}$. 
Since $S_a=f^{-1}(a)\cap Y$ is affine, 
there is a regular function $r$ on $S_a$ such that
$r(y_1)\neq r(y_2)$.

From the short exact sequence
$$ 0\longrightarrow 
 {\mathcal{O}}_Y
\longrightarrow 
 {\mathcal{O}}_Y
\longrightarrow 
 {\mathcal{O}}_{S_a}
\longrightarrow 
0,
$$
where the first map is defined by $f-a$, 
we
have surjective  map from  $H^0(Y, {\mathcal{O}}_Y)$ to 
$H^0(Z, {\mathcal{O}}_Z)$ since  $H^1(Y,  {\mathcal{O}}_Y)=0$.
 Lift $r$ to a regular function 
$R$ on $Y$, we find a function separating $y_1$ and $y_2$.
So $g$ is injective.   \\

\begin{remark} We   may  use the following alternating  construction of 
the morphism $g$
from $Y$         to   Spec$\Gamma(Y, {\mathcal{O}}_Y)$ 
as follows.

Since  $\kappa(D, X)=3$, 
we have a dominant morphism 
$h$  from $Y$  to
 $\Bbb{A}^d_{\Bbb{C}}$ 
defined by   three   algebraically independent nonconstant functions
on $Y$. Let $Z$ be the normalization of $\Bbb{A}^3$ in $Y$, then we
have a morphism $g$ from  $Y$  to $Z=$Spec$\Gamma(Y, {\mathcal{O}}_Y)$
and  $g$ is birational since $Y$ and $Z$ have  the same function field.
\end{remark}

 {\bf Step 6}.  $Y$ is affine.

$Proof.$ 
The above morphism  $g$  is birational since 
$Y$ and Spec$\Gamma(Y, {\mathcal{O}}_Y)$  have the same function field. 
 By Zariski's Main
Theorem  [27, Chapter 3, Section 9], $g$ is an open immersion from  
$Y$  
to ${\mbox{Spec}}\Gamma(Y, {\mathcal{O}}_Y)$.
By Neeman's theorem \cite{29},  $Y$ is affine since
$H^i(Y,  {\mathcal{O}}_Y)=0$ for all $i>0$. In fact, 
as a scheme,   $Y=$Spec$\Gamma(Y, {\mathcal{O}}_Y)$
since  they have the same coordinate ring.
\begin{flushright}
 Q.E.D. 
\end{flushright}

{\bf Proof of Theorem 1.4.} If $Y$ is an affine variety with dimension $d$, 
then it satisfies the following  four conditions  [15, Chapter 2, 
Section 1, 3, 6]:

(1) $Y$ satisfied condition A;

(2) $\kappa(D, X)\geq  1$;

(3) $ H^i(Y, {\mathcal{O}}_Y )=0 $ for all $i>0$. 

 Suppose now that the    quasi-projective variety 
 $Y$ satisfies the above   three conditions, we 
will prove that it is affine. We will use the induction on 
the dimension of $Y$. If $Y$ is a curve, it is trivial  [15,
Chapter 2, Section 4, Proposition 4.1]. 
If $Y$ is a surface, the proof is the same as
 in the proof of Theorem 1.1, Step 2.  
Assume  that the claim holds for every $d-1$ dimensional variety. 
Let dim$Y=d$. We may assume that $Y$ is normal   [14, Page 128,
Page 222].

{\bf Step 1}.  $\kappa (D, X)=d$.

$Proof.$   Since $\kappa (D, X)\geq 1$, we can construct a fibre space 
$f$: $X\rightarrow \bar{C}$, where $\bar{C}$  is a smooth projective curve. 
We also have the 
following commutative diagram  
 \[
  \begin{array}{ccc}
    Y                           &
     {\hookrightarrow} &
    X                                 \\
    \Big\downarrow\vcenter{%
        \rlap{$\scriptstyle{f|_Y}$}}              &  &
    \Big\downarrow\vcenter{%
       \rlap{$\scriptstyle{f}$}}      \\
C        & \hookrightarrow &
{\bar{C}},
\end{array}
\]
where  $C=f(Y)$   is   a smooth affine curve, 
$f$ is proper and surjective, 
every fibre of $f$ over $C$
 is connected.

 We will compute the $D$-dimension of  $X$ by Fujita's formula. 
For  the condition in  Fujita's  formula, 
see the proof of  Theorem 1.6.

Let  $\pi:  X'\rightarrow X $ be a proper
surjective morphism such  that $X'$ is   a   smooth projective variety
and every fibre is connected. In fact, 
$\pi$ is birational and generically finite.  Then  
we have [39, Chapter II, Theorem   5.13]
$$ \kappa(\pi ^*D, X')=\kappa(D, X).
$$
Consider a general  fibre $Y_t=f^{-1}(t)\cap Y$, 
$t\in C\subset\bar{C}$. 
By condition A, $Y_t$ is affine. Let 
$X_t=f^{-1}(t)$ and  $D_t$ be an effective boundary divisor
with support   $X_t-Y_t$, then
$\kappa(D_t, X_t)=d-1$.   
Pull this divisor  $D_t$ back to  $X'_t=\pi^{-1}(X_t) $, we have 
[39, Chapter II, Theorem 5.13]
$$  \kappa(\pi^*D_t, X'_t)=\kappa(D_t, X_t)=d-1. 
$$
Consider the new fibre space $X'\rightarrow  \bar{C}$, 
by Fujita's formula
we have 
$$ \kappa(\pi^*D, X')=\kappa(\pi^*D_t, X'_t)+1=d. 
$$ 
Therefore   $\kappa(D, X)=d$.

{\bf Step 2}. For any point $y$ in $Y$, there is a
 principal divisor $Z$ passing through $y$.

$Proof$.  We  may assume  $d\geq 2$.  Consider Ueno's   fibre space  [39, Page 46]
in Step 4, proof of Theorem 1.1
\[
  \begin{array}{ccc}
    Y                           &
     {\hookrightarrow} &
    X                                 \\
    \Big\downarrow\vcenter{%
        \rlap{$\scriptstyle{f|_Y}$}}              &  &
    \Big\downarrow\vcenter{%
       \rlap{$\scriptstyle{f}$}}      \\
C        & \hookrightarrow &
{\Bbb{P}}^1.
\end{array}
\]
We  know that  every  fibre  of  $f$  is  connected. 
Let $f(y)=a$, then $a\in f^{-1}(a)=X_t$. Let  $Z=Y_t=X_t\cap Y$, 
then  $Z=\{f-a=0, f\in H^0(Y, {\mathcal{O}}_Y)\}$
is a principal divisor passing through $y$.

      By condition  A, $Z$  is affine.\\

{\bf Step 3}. $Y$  contains no complete curves.

$Proof$.   If  $Y$  has  a  complete curve  $F$, we may assume that 
it is connected. Let  $p\in  F$  be a point on  $F$.  Let $f(p)=a$,
where $f$ is the morphism in Step 1. Then  $a$ is a point on $C=f(Y)$.
Let $X_a=f^{-1}(a)$, $Y_a=X_a\cap  Y$, then  $p$  is a point on the 
open fibre $Y_a$. Since  $F$ is connected  and $X_a$ is also connected,
 $F\subset  Y_a$. But $Y_a$  is a closed subvariety of  $Y$, so
$Y_a$ is affine by condition   A. This is a contradiction
since  an affine variety does not have  any complete curve.
 Thus  $Y$
contains  no complete curves. \\

{\bf Step 4}.  There is an injective  morphism $g$ from $Y$ to 
${\mbox{Spec}}\Gamma(Y, {\mathcal{O}}_Y)$.

$Proof$. We again define the morphism from 
$Y$ to 
${\mbox{Spec}}\Gamma(Y, {\mathcal{O}}_Y)$
as in the Step 5, proof of Theorem 1.1.
We send a closed irreducible subvariety on $Y$ (as a quasi-projective
variety) to its defining prime ideal. We know this morphism 
$g$ is injective if and only if for any two distinct points
$y_1$ and $y_2$ on $Y$, there is a regular function $H$
on $Y$ such that  $H(y_1)\neq H(y_2)$. 
Let $Z$ be the  principal divisor in Step 2 passing though 
$y_1$ and $f\in \Gamma(Y, {\mathcal{O}}_Y)$ be its defining
function. If $f(y_2)\neq f(y_1)$, then we are done. Suppose 
 $f(y_2)=f(y_1)=0$. From the exact sequence
$$ 0\longrightarrow 
 {\mathcal{O}}_Y
\longrightarrow 
 {\mathcal{O}}_Y
\longrightarrow 
 {\mathcal{O}}_Z
\longrightarrow 
0,
$$
we can lift any regular function from $Z$ to $Y$. 
In particular, we can choose $h\in \Gamma(Z, {\mathcal{O}}_Z)$
such that $h(y_1)=0$ but $h(y_2)=1$. Thus there is a regular function $H$
on $Y$ such that $H(y_1)\neq H(y_2)$.  
\\

{\bf Step 5}. $Y$ is affine. 

$Proof$. By Zariski's Main Theorem, 
$Y$ is isomorphic to an open subset of 
Spec$\Gamma(Y, {\mathcal{O}}_Y)$.
By Neeman's theorem, $Y$ is affine. 
\begin{flushright}
 Q.E.D. 
\end{flushright}

{\bf Proof of Theorem 1.5.} If $Y$ is affine, then 
$H^i(Y, \Omega^j_Y)=0$   for all $j\geq 0$,  $i>0$
and $\kappa(D, X)=n$  [15, Chapter 2, Theorem 1.1, Theorem 6.1].
Let $Y$ be a smooth quasi-projective variety with 
$H^i(Y, \Omega^j_Y)=0$   for all $j\geq 0$,  $i>0$
and $\kappa(D, X)=d$, we  need to  prove that it is affine. 
We will use the induction on the dimension of $Y$.

The claim holds for
curves [15, Chapter 2, proposition 4.1] and surfaces
by the proof of Theorem 1.1, Step 2. Assume that the claim holds for
$d-1>0$ dimensional smooth varieties. Let dim$Y=d$.

{\bf Step 1}.
We first claim that for any point $y$ in $Y$, 
there is a smooth prime principal 
divisor $Z=\{f=0, f\in H^0(Y, {\mathcal{O}}_Y) \}$ passing through $y$.

$Proof.$
Since $\kappa(D, X)=$dim$X=d>0$, there is $a>0$
such  that for all  $m\gg 0$, $h^0(X, {\mathcal{O}}_X(mD))\geq a m^d$
[39, Theorem 8.1]. 
Recall our notation  in Section 2.
Choose a basis $\{f_0, f_1, \cdot \cdot\cdot, f_n\}$
     of the linear space 
     $H^0(X, {\mathcal{O}}_X(mD))$, it defines a rational 
     map 
     $\Phi _{|mD|}$
     from $X$ to the projective space 
     ${\Bbb{P}}^n$ by sending a point $x$ on $X$ to
     $(f_0(x), f_1(x), \cdot \cdot\cdot, f_n(x))$ in ${\Bbb{P}}^n$. 
By definition of $D$-dimension,
$$ \kappa (D, X)= \max_m\{\dim (\Phi _{|mD|}(X))\}=d. 
      $$
Let ${\Bbb{C}}(X)$ be the function field of $X$. Let 
$$  R(X, D)=\oplus_{\gamma\geq 0}H^0(X, {\mathcal{O}}_X(\gamma D)) 
$$
be the graded $\Bbb{C}$-domain  and $R^*\subset R$
the  multiplicative  subset  of all nonzero
homogeneous  elements. Then the quotient ring $R^{*-1}R$  
is a graded  $\Bbb{C}$-domain  and its degree 0 part 
$(R^{*-1}R)_0$  is a  field  denoted  by  $Q((R))$, or $Q((X, D))$, i.e., 
$$ Q((R))=Q((X, D))=(R^{*-1}R)_0.
$$ 
It is easy to see that  $Q((X, D))=Q((X, aD))$ 
and $\kappa(D, X)=\kappa(aD,  X)$ for all
positive  integer $a$. 
Notice that the quotient field $Q(R)$ of $R$ is different from 
$Q((R))$. 

By [26, Proposition 1.4], 
the function field of 
the image $\Phi _{|mD|}(X)$ is  $Q((X, D))$, i.e., 
${\Bbb{C}}(\Phi _{|mD|}(X))=Q((X, D))$. 
By [26, Proposition 1.9], 
$\Phi _{|mD|}$ is birational for all $m\gg 0$
since  $D$ is effective and    $\kappa(X, D)=$dim$X$. 
So the two function fields 
${\Bbb{C}}(\Phi _{|mD|}(X))$   and  
${\Bbb{C}}(X)$
are isomorphic
$$ {\Bbb{C}}(\Phi _{|mD|}(X))= {\Bbb{C}}(X). 
$$ 
Therefore  
$$   {\Bbb{C}}(X)=Q((X, D)).
$$
This implies that any rational function $r$
on $X$ can be written as a quotient $f/g$
of two regular functions $f$ and $g$ on $Y$ with $g\neq 0$, 
where there is some $m>0$   such  that 
 $f, g\in  H^0(X, {\mathcal{O}}_X(mD))$. 
If $r$ is regular on an open subset $U$
in $Y$, then $r=0$ on a closed subset $A$ of $U$
gives $f=0$  on $A$. 

For any point $y$ on $Y$,  there is a smooth prime divisor
$Z$  passing through  $y$. Since $Y$ is smooth, it is 
locally factorial.  Let $r$ be the local defining 
function of $Z$ near $y$, then there is a regular function
$f$ on $Y$ such that $f$ also defines $Z$ near $y$. We  define  the smooth
prime divisor $Z=\{f=0\}.$

{\bf Step 2}. The above  smooth principal divisor $Z$
satisfies  the same vanishing  condition $H^i(Z, \Omega_Z^j)=0$
for all $i>0$ and $j\geq 0$. Moreover, 
$Z$ is affine. 

$Proof$.  Let $f$ be the defining function of 
$Z$, then $$ 0\longrightarrow 
\Omega^j_Y
\longrightarrow 
\Omega^j_Y
\longrightarrow 
\Omega^j_Y|_Z
\longrightarrow 
0, 
$$
where the first map is defined by $f$. 
So $H^i(Z, \Omega^j_Y|_Z)=0$  since 
$H^i(Y, \Omega^j_Y)=0$   for all $j\geq 0$,  $i>0$. In particular, 
$H^i(Z, {\mathcal{O}}_Z)=0$ for all $i>0$. 
Since both $Y$ and $Z$ are smooth, 
We have an exact sequence  [14, Chapter 2, Section 8; 
12, Page 157],
$$0\longrightarrow 
\Omega^{j-1}_Z
\longrightarrow 
\Omega^j_Y|_Z
\longrightarrow 
\Omega^j_Z
\longrightarrow 
0. 
$$ 
So  $H^i(Z, \Omega^j_Z)=0$   for all $j\geq 0$,  $i>0$.

By Theorem 5.11 and Theorem 5.13, \cite{39}  
$$  d= \kappa(D, X)\leq  \kappa(D_t, X_t )+1,
$$ 
where $X_t=f^{-1}(t)$ is a general fibre and $D_t$ is the boundary divisor 
on $X_t$ supported in $X_t-Y$ for $t\in C \subset \Bbb{C}$. 
So $\kappa(D_t, X_t)=d-1$. 
Let $\bar{Z}=X_0$ be the  completion of $Z$  in $X$,
let $D_0$ be the effective boundary divisor supported in $X_0-Y$,
then  by Grauert's  upper semi-continuity theorem
[39, Chapter 1], $\kappa(D_0, \bar{Z})=d-1$. 

By the inductive assumption, $Z$ is affine.\\

{\bf Step 3}.   $Y$ is affine.

$Proof$.  
We construct the  morphism
$g$
from $Y$ to Spec$\Gamma(Y, {\mathcal{O}}_Y)$ as  in Step 5, proof of 
Theorem 1.1. $g$ is birational since  $Y$ and 
 Spec$\Gamma(Y, {\mathcal{O}}_Y)$ 
have the same   function field. It is also injective. In fact, 
for any two distinct points $y_1$ and $y_2$ on $Y$ there is
a regular function $R$ on $Y$ such that $R(y_1)\neq R(y_2)$.
This is  the consequence of the affineness
of the smooth principal divisor 
$Z=\{f=0\}$
passing through $y_1$.  More precisely, 
if  $y_2\notin Z$, then  $f(y_1)\neq f(y_2)$, we are done.
Assume  $y_2\in Z$, i.e., $f(y_1)=f(y_2)=0$. 
By the affineness of $Z$ and the short exact sequence
$$ 0\longrightarrow 
 {\mathcal{O}}_Y
\longrightarrow 
 {\mathcal{O}}_Y
\longrightarrow 
 {\mathcal{O}}_Z
\longrightarrow 
0,
$$
we can lift a regular function $r$ from $Z$ to $Y$,
where $r(y_1)\neq r(y_2)$. 
So  if $R$ is the lifting, then 
$R(y_1)\neq R(y_2)$. This proves the injectivity of    
$g$.  By Zariski's Main Theorem and Neeman's theorem, $Y$ is affine.
\begin{flushright}
 Q.E.D. 
\end{flushright}
{\bf Proof of Theorem 1.6.}  We will prove that $\kappa(D, X)=$dim$X=n$ then 
 apply
Theorem 1.5. 

Let $f: Y\longrightarrow  C$ be the surjective morphism 
from $Y$ to the smooth affine curve $C$. Then $f$ gives a 
rational map from $X$ to  the completion   $\bar{C}$
of $C$.   
Resolve the
indeterminacy of $f$ on the boundary $X-Y$.
We may replace $X$ by its suitable blowing up and assume that
$f:X\rightarrow  \bar{C}$  is surjective and proper morphism. 
Notice that this  procedure  does not  change $Y$.
$Y$ is still an open subset of $X$. 
By Stein factorization, we may assume that every fibre is 
connected and general fibre is smooth.
Pick a point $t_1\in \bar{C}-C$, then 
$\kappa(t_1, \bar{C})=1$.  For a general point 
$t\in C$, by the Riemann-Roch formula, there is  a positive integer $m$,
 such that 
$h^0(\bar{C}, {\mathcal{O}}(mt_1-t))>1$.  Let $s$ be a nonconstant section of 
$H^0(\bar{C}, {\mathcal{O}}(mt_1-t))$, then
$$ {\mbox{div}} s +mt_1-t\geq 0.
$$
Pull it back to $X$, we have 
$$  f^*({\mbox{div}} s +mt_1-t)={\mbox{div}} f^*(s)+mf^*(t_1)-f^*(t)\geq 0.
$$
Let $D_1=f^*(t_1)$ and $F=f^*(t),$ then  
$h^0(X, {\mathcal{O}}_X(mD_1-F))>0$.  Choose an effective divisor 
$D$ with support $X-Y$ such that  $D_1\leq   D$, then we  have
$$ h^0(X, {\mathcal{O}}_X(mD-F))\geq  h^0(X, {\mathcal{O}}_X(mD_1-F))>0. 
$$
Since  $F|_Y$ is a smooth affine subvariety of codimension 1
 [15, Chapter 2, Proposition 4.1],
 $\kappa(D|_F, F)=n-1$.
By Fujita's equality, 
$$\kappa(D, X)=\kappa(mD, X)=\kappa(mD|_F, F)+\kappa(t_1, \bar{C})=n.$$
By Theorem  1.5, $Y$ is affine.
\begin{flushright}
 Q.E.D. 
\end{flushright}

It is easy to see that Theorem 1.7 follows from Theorem 1.6.

\begin{remark}  In Theorem 1.5, if we drop one of the two  conditions
in the theorem,  the  theorem is no longer true.

Example 1.  A threefold  $Y$  with  $H^i(Y, \Omega^j_Y)=0$
for all $i>0$  and $j\geq 0$  but  $\kappa(D, X)=1$.

 Let $C$ be an elliptic curve and $E$ the unique nonsplit 
      extension of $\mathcal{O}$$_C$ by itself.  
      Let ${Z=\Bbb{P}}_C(E)$ and  $D$ be the canonical section,
then  $H^i(S, \Omega^j_S)=0$ for all $i>0$ and $j\geq 0$, where  $S=Z-D$
\cite{24}.
Let  $F$ be a smooth affine curve and $Y=S\times F$, then  
$H^i(Y, \Omega^j_Y)=0$ by   K$\mbox{\"{u}}$nneth   formula \cite{35}. Let
$X$ be the closure of $Y$ and $D$ be the effective boundary 
divisor, then $\kappa(D, X)=1$ \cite{42}. So $Y$ is not affine. 

Example 2.  A threefold  $Y$  with  $\kappa(D, X)=3$
but  $Y$  is not affine. 

For instance, remove a hyperplane section $H$ and a line
$L$ from ${\Bbb{P}}^3$, where $L$ is not contained in $H$. 
Let $Y={\Bbb{P}}^3-H-L$.  Then $Y$ contains no complete curves.
 Let 
$f:X\rightarrow {\Bbb{P}}^3$
be the blowing up of ${\Bbb{P}}^3$ along $L$. Then $X$ is a  smooth 
projective threefold and  $Y$ is an open subset of $X$. 
Let $D=f^{-1}(H)+E$, where $E$ is the exceptional divisor.
Then by Iitaka's formula  in Section 2,
$\kappa (D, X)=\kappa (H, {\Bbb{P}}^3)=3$. But $Y$ is not affine 
since  the boundary ${\Bbb{P}}^3-Y$ is not of pure codimension
1  [15,  Chapter 2, Proposition 3.1].

\end{remark}

\end{document}